\documentclass[11pt]
{amsart}
\usepackage{amsthm,amsmath,amssymb,amsfonts}
\usepackage[backref,colorlinks]{hyperref}
\usepackage[margin=1.4in]{geometry}

\newtheorem{theorem}{Theorem}
\newtheorem{thm}[theorem] {Theorem}
\newtheorem{lemma}[theorem]{Lemma}
\newtheorem{lem}[theorem]{Lemma}

\newtheorem{cor}[theorem]{Corollary}

\newtheorem{fact}[theorem]{Fact}
\newtheorem{conj}[theorem]{Conjecture}

\newtheorem{obs}[theorem]{Observation}
\newtheorem{defi}[theorem]{Definition}
\newtheorem{question}[theorem]{Question}

\let\eps=\varepsilon
\let\theta=\vartheta
\let\rho=\varrho
\let\sigma=\varsigma

\let\polishlcross=\l
\def\l{\ifmmode\ell\else\polishlcross\fi}

\newcommand{\diam}{\text{diam}}

\def\cM{{\mathcal M}}

\def\cF{{\mathcal F}}

\def\a{\alpha}

\def\d{\delta}

\def\l{\lambda}

\def\r{\rho}

\def\s{\sigma}

\def\Pr{\mbox{{\bf Pr}}}

\def\Hcs{Hamilton cycles}

\newcommand{\proofstart}{{\bf Proof\hspace{2em}}}

\newcommand{\proofend}{\hspace*{\fill}\mbox{$\Box$}\vspace*{7mm}}

\newcommand{\Exp}{\mbox{\bf E}}

\newcommand{\rdown}[1]{{\mbox{$ \lfloor #1 \rfloor $}}}

\newcommand{\rt}{\right}

\newcommand{\lt}{\left}

\def\cH{{\mathcal H}}

\begin{document}

\title[On covering expander graphs by Hamilton cycles]{On covering expander graphs by Hamilton cycles}

\author[R.~Glebov]{Roman Glebov}
\address{Institut f\"{u}r Mathematik, Freie Universit\"at Berlin, Arnimallee 3-5, D-14195 Berlin, Germany}
\email{glebov@mi.fu-berlin.de}
\thanks{The first author was supported by DFG within the research training group "Methods for Discrete Structures".}
\author[M.~Krivelevich]{Michael Krivelevich}
\address{School of Mathematical Sciences,
Sackler Faculty of Exact Sciences,
Tel Aviv University,
Tel Aviv 69978,
Israel}
\email{krivelev@post.tau.ac.il}
\thanks{The second author was supported in part by a USA-Israel BSF grant and by a grant from the Israel Science
Foundation.}
\author[T.~Szab\'o]{Tibor Szab\'o}
\address{Institut f\"{u}r Mathematik, Freie Universit\"at Berlin, Arnimallee 3-5, D-14195 Berlin, Germany}
\email{szabo@mi.fu-berlin.de}

\date{\today}
\begin{abstract}
The problem of packing Hamilton cycles in random and pseudorandom
graphs has been studied extensively. In this paper, we look at the
dual question of covering all edges of a graph by Hamilton cycles
and prove that if a graph with maximum degree $\Delta$ satisfies
some basic expansion properties and contains a family of
$(1-o(1))\Delta/2$ edge disjoint Hamilton cycles, then there also
exists a covering of its edges by $(1+o(1))\Delta/2$ Hamilton
cycles. This implies that for every $\a
>0$ and every $p \geq n^{\a-1}$ there exists a covering of all edges
of $G(n,p)$ by $(1+o(1))np/2$ Hamilton cycles asymptotically almost
surely, which is nearly optimal.
\end{abstract}

\maketitle

\setcounter{footnote}{0}
\renewcommand{\thefootnote}{\fnsymbol{footnote}}

\section{Introduction}

For an $r$-uniform hypergraph $G$ and a family $\cF$ of its subgraphs,
we call a family $\cF'\subset \cF$ an $\cF$-{\em decomposition} of $G$
if every edge of $G$ is contained in {\em exactly} one of the hypergraphs from $\cF'$.
We call a family $\cF'\subset \cF$  an $\cF$-{\em packing} of $G$, if every edge of $G$ is contained in {\em at most} one of the hypergraphs from $\cF'$.
Naturally, one tries to {\em maximize} the size of an $\cF$-packing of $G$.
The dual concept is that of an $\cF$-covering:
a family $\cF'$ is called an $\cF$-{\em covering} of $G$, if every edge of $G$ is contained in {\em at least} one of the hypergraphs from $\cF'$.
Here the {\em minimum} size of an $\cF$-covering of $G$ is sought.

Decompositions, packings and coverings are in the core of combinatorial research (see~\cite{Furedi} for a survey).
One of the most famous problems in this area was the conjecture of Erd\H{o}s and Hanani~\cite{EH},
dealing with the case when $G$ is the complete $r$-uniform hypergraph on $n$ vertices and
$\cF$ is the family of all $k$-cliques in $G$ for some $k\geq r$.
Clearly, if there was an $\cF$-decomposition of $G$, its size would be $\binom{n}{r}/\binom{k}{r}$.
Hence $\binom{n}{r}/\binom{k}{r}$ is an upper bound on the size of a largest $\cF$-packing of $G$
and a lower bound on the size of a smallest $\cF$-covering of $G$.
Erd\H{o}s and Hanani conjectured both inequalities to be asymptotically tight for constant $r$ and $k$,
i.e., that a largest $\cF$-packing of $G$ and a smallest $\cF$-covering of $G$ are asymptotically equal
to each other.
R\"odl~\cite{Rodl} verified the conjecture by one of the first applications of the nibble method.
Observe that in this setting, the two parts of the conjecture are trivially equivalent.
The reason for this is that the size of the elements of the family $\cF$ does not grow with $n$:
from a packing $\cF'$ of $G$ of size $(1-\eps)\binom{n}{r}/\binom{k}{r}$,
one obtains a covering $\cF''$ of $G$ of size $(1+\eps\binom{k}{r})\binom{n}{r}/\binom{k}{r}$
by simply taking additionally one $k$-clique for every $r$-edge that was not contained in any clique from $\cF'$.

In this paper we study a covering problem where the sets in $\cF$ grow with $n$ and the above equivalence is
not entirely clear.  Let $r=2$, so our objects are usual graphs.
For a graph $G$, the considered family $\cH=\cH (G)$ is the family of all {\Hcs} of $G$.
The corresponding concepts of decomposition, packing, and covering are called
{\em Hamilton decomposition}, {\em Hamilton packing} and {\em Hamilton covering}, respectively.

The most well-known fact about Hamilton decompositions is the nearly folklore result of Walecki
(see e.g.~\cite{Walecki}),
stating that for every odd $n$, the complete graph $K_n$ has a Hamilton decomposition.
In general, however, not many graphs are known to have a Hamilton decomposition; the interested reader is referred to~\cite{cyclesandrays}.

Given that the minimum degree of a Hamilton cycle is $2$, the
maximum size of a Hamilton packing of a graph with minimum degree
$\d$ is $\lfloor\delta/2\rfloor$. Interestingly, the random graph
$G(n,p)$ seems to match this bound tightly. There has been an
extensive research considering Hamilton packings of the random graph
$G(n,p)$. A classic result of Bollob\'as~\cite{bollobas} and
Koml\'os and Szemer\'edi \cite {KS} states that as soon as the
minimum degree of the random graph is $2$, it contains a Hamilton
cycle a.a.s. This result was extended by Bollob\'as and
Frieze~\cite{BF} who showed that we can replace $2$ by $2k$ for any
constant $k$ and obtain a Hamilton packing of size $k$ a.a.s. Frieze
and Krivelevich~\cite{FK1} proved that for every constant and
slightly subconstant $p$, $G\sim G(n,p)$ contains a packing of
$(1+o(1))\delta(G)/2$ Hamilton cycles a.a.s. They also conjectured
that for every $p=p(n)$ there exists a Hamilton packing of $G\sim
G(n,p)$ of size $\lfloor\delta(G) /2\rfloor$ (and, in case
$\delta(G)$ is odd, an additional (disjoint to the Hamilton cycles
of the packing) matching of size $\lfloor n/2\rfloor$) a.a.s. Frieze
and Krivelevich~\cite{FK2} proved their conjecture as long as $p=
(1+o(1))\ln n/n$, which was extended to the range of $p\leq1.02\ln
n/n$ by Ben-Shimon, Krivelevich and Sudakov~\cite{BSKS}. Meanwhile,
Knox, K\"uhn and Osthus~\cite{KKO1} extended the result
from~\cite{FK1} to the range of $p=\omega(\ln n /n)$, and then
proved the conjecture for
$\ln^{50}n/n<p<1-n^{-1/4}\ln^9n$~\cite{KKO2}. Very recently it was
also proven by Krivelevich and  Samotij~\cite{KrSa}  that there
exists a positive constant $\eps>0$ such that for the range of $\ln
n/n\leq p \leq n^{\eps-1}$, $G\sim G(n,p)$ contains a Hamilton
packing of size $\lfloor\delta(G) /2\rfloor$ a.a.s., implying the
conjecture in this range of $p$ up to the existence of the
additional matching.

To the best of our knowledge the dual concept of Hamilton covering of $G(n,p)$ has not been studied.
Obviously, the size of any Hamilton covering of graph $G$ is at least $\lceil\Delta(G)/2\rceil$,
where $\Delta (G)$ denotes the maximum degree of $G$.
Recall that for $p=p(n)\gg \log n/n$, $\Delta(G(n,p))=(1+o(1))np = \delta( G(n,p))$ a.a.s.,
hence the minimum size of a Hamilton cover and the maximum size of a Hamilton packing have
a chance to be asymptotically equal.
We prove that this, in fact, is the case for the range $p > n^{\alpha -1}$ where $\alpha >0$ is an arbitrary small
constant.

\begin{thm}
\label{corrandom} For any $\alpha > 0$, for $p \geq n^{\alpha-1}$
a.a.s. $G(n,p)$ can be covered by $(1+o(1))np/2$ Hamilton cycles.
\end{thm}


\subsection{Pseudorandom setting}

Our argument will proceed in an appropriately chosen pseudorandom setting.

By the {\em neighborhood} $N(A)$ of a set $A$, we mean all the vertices outside $A$ having at least one neighbor in $A$.
Note that we explicitly exclude $A$ from $N(A)$.

The following definition contains the most important notions of the paper.
\begin{defi}
We say that a graph $G$ has the {\em small expander property} $S(s,g)$ with {\em expansion factor} $s$ and {\em boundary} $g$,
if for any set $A\subset V(G)$ of size $|A|\leq g$, the neighborhood of $A$ satisfies $|N(A)|\geq |A|s$.

We say that $G$ has the {\em large expander property}  $L(l)$ with {\em frame} $l$, if there is an edge between any two disjoint sets $A,B\subset V(G)$ of size $|A|,|B|\geq l$.

We call a graph $(s,g,l)$-{\em expander}, if it satisfies properties $S(s,g)$ and $L(l)$.

We refer to an $\lt(s,  \frac{4n\ln s}{s\ln n}, \frac{n\ln s}{3000\ln n}\rt)$-expander on $n$ vertices briefly as an $s${\em -expander}.
\end{defi}

Notice that the expander-property is monotone in all three parameters,
meaning that  every $(s,g,l)$-expander is also an $(s-1, g,l)$-, an $(s,g-1,l)$- and an $(s,g,l+1)$-expander.

The proof of Theorem~\ref{corrandom} is based on the following result.

\begin{thm}
\label{teo}
For every constant $\alpha>0$ and $n$ large enough,
every $n^{\alpha}$-expander graph $G$ on $n$ vertices with a Hamilton packing of size $h$ has a
Hamilton covering of size at most $h+\frac{28000(\Delta(G)-2h)}{\alpha^4}$.
\end{thm}


\subsection{Structure of the paper and outline of the proofs}

In Section~\ref{pf1}, we prove Theorem~\ref{teo} in the following main steps:
\begin{itemize}
\item
in Lemma~\ref{inducedexpander} we show that the small expander property is ``robust'' in the sense that after deleting a small linear-size set of vertices from a graph satisfying $S$, we still have a large subgraph satisfying $S$ with slightly worse parameters;
\item
in Lemma~\ref{largematching} we show that vertex-disjoint paths in an expander can be concatenated without loosing or gaining too many edges;
\item
in Lemma~\ref{St}, Fact~\ref{good} and Lemma~\ref{path+}, we learn how to apply the rotation-extension technique, developed by P\'osa~\cite{Posa}, without loosing too many important edges;
\item
Lemma~\ref{mainlemma} contains the main proof of the paper.
There we show applying the previous technical statements that a small matching in an expander can mostly be covered by a Hamilton cycle of the same graph.
\item
After having digested the statement of Lemma~\ref{mainlemma} in Corollary~\ref{oldteo}, we prove Theorem~\ref{teo} in $6$ easy lines.
\end{itemize}

Section~\ref{pf2} contains the proof of Theorem~\ref{corrandom}.
There, we first prove in Lemma~\ref{expander} that $G(n,p)$ is a $\sqrt[5]{np}$-expander a.a.s.,
and prove Theorem~\ref{corrandom} using this fact and the result from~\cite{KKO1}.

And finally, in Section~\ref{concl}, we give some concluding remarks leading to further open questions.

In general, we may drop floor and ceiling signs to improve the readability when they do not influence the asymptotic statements.

\section{Proof of Theorem~\ref{teo}}\label{pf1}

The maximum over all pairs of vertices $x,y\in V(G)$ of the length of a shortest $xy$-path in a graph $G$
is called the {\em diameter} of $G$ and is denoted by $\diam (G)$.
We start with an observation showing that graphs with appropriate expander properties have small diameter.
\begin{obs}
\label{diam}
Any $n$-vertex graph satisfying $S(s,g)$ and $L(l)$ for some $s,g,l$ with $s>1$ and $l\leq sg$ has diameter at most $2\ln n/\ln s + 3$.
\end{obs}

\proofstart
Let $G$ be a graph on $n$ vertices.
By the small expansion property of $G$, we know for every vertex $x\in V(G)$ that $|N(x)|\geq s$.
In the following, for any $x\in V(G)$ and any $i\in \mathbb{N}$,
let us denote by $B_i(x)$ the set of all vertices with distance at most $i$ to $x$,
i.e., $B_i(x)$ contains all those vertices $y\in V(G)$, for which there exists an $xy$-path of length at most $i$.
Inductively, as long as for an $i\in \mathbb{N}$ and an $x\in V(G)$ it holds that $|B_i(x)|\leq g$, we obtain $|B_{i+1}(x)|\geq s^{i+1}$.
Now, the last index $i$ such that $|B_i(x)|\leq g$
is obviously at most $\left\lfloor \ln n/\ln s\right\rfloor -1$, since $sg\leq n$.
Hence, $|B_{\left\lfloor\ln n/\ln s\right\rfloor}(x)|\geq g$,
and $|B_{\left\lfloor\ln n/\ln s\right\rfloor+1}(x)|\geq sg\geq l$.
Thus, for any two vertices $x,y\in V(G)$ we know that both sets $B_{\left\lfloor\ln n/\ln s\right\rfloor+1}(x)$
and $B_{\left\lfloor\ln n/\ln s\right\rfloor+1}(y)$ have at least $l$ vertices.
This guarantees that by property $L(l)$ of $G$, either these sets are not disjoint or
there exists an edge between these sets.
Each of these facts implies an $x,y$-walk of length at most $2\left\lfloor\ln n/\ln s\right\rfloor+3$,
and the observation follows.
\proofend

The following lemma shows that if we have a graph satisfying the small expansion property
and we remove an arbitrary subset of small size from the vertex set,
then with the additional removal of an even smaller subset
we can recover some of the small expander property again.

\begin{lem}[Induced expander lemma]
\label{inducedexpander}
For any $s,g$, every graph $G=(V,E)$ satisfying $S(s,g)$ has the following induced expander property.
For every $D\subset V$ of size $|D|\leq \frac{gs}{4}$, there exists a set $Z\subset V$ of size $|Z|\leq \frac{2|D|}{s}$,
such that the graph $G[V\setminus(D\cup Z)]$ satisfies $S(\frac{s}{2}, \frac{g}{2})$.
\end{lem}

\proofstart
Let $Z$ be a largest set in $V\setminus D$ among subsets of $V\setminus D$ of size at most $g$
not satisfying the small expander property with expansion factor $s/2$ in $G[V\setminus D]$,
meaning that  $|N(Z)\setminus D|< |Z|\frac{s}{2}$
(assuming there exists such a set; otherwise we are done by setting $Z=\emptyset$ in the statement of the lemma).
We denote $U=V\setminus(D\cup Z)$ and remember that
\[|N(Z)\cap U|<  \frac{|Z|s}{2}.\]
Thus, by the property $S(s,g)$ of $G$,
\[|D|+|N(Z)\cap U|\geq |N(Z)|\geq |Z|s,\]
implying that
\[|Z|\leq \frac{2|D|}{s} \leq \frac{g}{2}.\]

Assume now for the sake of contradiction that $G[U]$ does not satisfy $S(\frac{s}{2}, \frac{g}{2})$,
i.e. there exists an $A\subset U$ with $|A|\leq \frac{g}{2}$ and $|N(A)\cap U| < |A|\frac{s}{2}$.
Then the set $A\cup Z\subset V\setminus D$ satisfies both properties
\[|A\cup Z|\leq g\]
and
\[|N(A\cup Z)\cap U|\leq |N(A)\cap U|+|N(Z)\cap U| < \frac{|A\cup Z|s}{2},\]
contradicting the assumption that $Z$ is a largest set in $V\setminus D$ with these properties.
\proofend


%
%
%
%
%
%
%
%
%
%
%
%

The following concept allows us to concatenate paths in an appropriate way.

By the $k${\em -end} of a path we mean the at most $2k$ vertices of this path's vertex set with distance at most $k-1$ to one of the two endpoints,
whereas the endpoints have distance $0$ to themselves and so are part of any $k$-end with $k\geq 1$.
We call a path {\em non-trivial} if uts length is at least one.
Given a family $\cM$ of non-trivial vertex-disjoint paths,
we call a family $\cM'$ of non-trivial vertex-disjoint paths a $(d,k)${\em -extension of $\cM$},
if we have that
\begin{itemize}
\item $\mu:=|\cM|-|\cM'|\geq 0$,
\item $|\bigcup_{P\in \cM} E(P)\setminus \bigcup_{P'\in \cM'} E(P')|\leq 2(k-1)\mu$,
and
\item $|\bigcup_{P'\in \cM'} E(P')\setminus \bigcup_{P\in \cM} E(P)|\leq (d+2)\mu$.
\end{itemize}
Informally speaking, we neither gain nor loose too many edges relative to the decrease in the number of paths,
while passing from $\cM$ to $\cM$.

To clarify the notation, when speaking about the {\em size} $|\cM|$ of a family of non-trivial
vertex-disjoint paths $\cM$, we mean the number of these paths.
Note that any extension of $\cM$ has at most as many paths as $\cM$.
Notice that the definition of extension is transitive,
i.e., if for some pair $(d,k)$, $\cM'$ is a $(d,k)$-extension of $\cM$
and $\cM''$ is a $(d,k)$-extension of $\cM'$, then $\cM''$ is a $(d,k)$-extension of $\cM$.
Furthermore, the relation is also reflexive, i.e., $\cM$ is a $(d,k)$-extension of itself for any $k\geq 1$ and every $d\geq 0$.
We say that $\cM'$ is a size-minimum $(d,k)$-extension of $\cM$,
if every $(d,k)$-extension $\cM''$ of $\cM$ has size at least $|\cM''|\geq |\cM'|$.
Notice that the transitivity and reflexivity of extensions imply that then $\cM'$ is a size-minimum $(d,k)$-extension of itself.

We use the notation $V(\cM):=\bigcup _{P\in \cM} V(P)$ for the set of all vertices appearing in one of the paths of $\cM$
and $E(\cM):=\bigcup _{P\in \cM} E(P)$ for the corresponding edge set.

For our applications we will mostly be interested in size-minimum extensions.
The following lemma provides us with two basic properties of size-minimum extensions:
namely that they do not contain very short paths and that there are no
short paths between the ends of two distinct paths from such an extension.

\begin{lemma}
\label{extension}
Let $G$ be a graph, $k\geq 1$ and $d\geq 0 $ arbitrary integers and $ \cM$ a family of non-trivial vertex-disjoint paths, which
is a size-minimum $(d,k)$-extension of itself.
Then the following holds:
\begin{enumerate}
\item
there exists no path in $ \cM$  of length less then $2k-1$ and
\item
for every two distinct paths $P_1, P_2\in \cM$, for every vertex $x$ in the $k$-end of $P_1$ and every vertex $y$ in the $k$-end of $P_2$,
for every $a\in N(x)\setminus V(\cM)$ and $b\in N(y)\setminus V(\cM)$, there exists no $ab$-path of length at most $d$ in $G - V(\cM)$.
\end{enumerate}
\end{lemma}

\proofstart
1. Suppose for the sake of contradiction that there exists a path $P$ of length at most $2k-2$ in $\cM$.
Delete this path from $\cM$ and call the resulting family of non-empty vertex-disjoint paths $\cM'$.
Then since $|\cM|-|\cM'|=1$, we obtain
\[|E(\cM)\setminus E(\cM')|=|E(P)|\leq 2k-2=2(k-1)\lt(|\cM|-|\cM'|\rt),\]
and
\[|E(\cM')\setminus E(\cM)|=0\leq (d+2)\lt(|\cM|-|\cM'|\rt).\]
Hence, $\cM'$ is a $(d,k)$-extension of $\cM$,
contradicting the minimality of $\cM$ as a $(d,k)$-extension of itself.

2. Suppose to the contrary that there exist two distinct paths $P_1, P_2\in \cM$
with a vertex $x$ in the $k$-end of $P_1$, a vertex $y$ in the $k$-end of $P_2$,
and vertices $a\in N(x)\setminus V(\cM)$ and $b\in N(y)\setminus V(\cM)$
such that in $G - V(\cM)$, there is an $ab$-path of length at most $d$.
Let us call this path $P_3$.

The vertex $x$ splits $P_1$ into two subpaths, the shorter one has length at most $k-1$.
Let us call the longer one $P_x$ and construct $P_y$ from $P_2$ analogously.
By connecting the paths $P_x$, $P_3$ and $P_y$ via the edges $xa$ and $by$,
we obtain a new path $P'$.
Replace the two paths $P_1$ and $P_2$ in $\cM$ by $P'$ and call the resulting family of non-trivial
vertex-disjoint paths $\cM'$.
Then $\mu=|\cM|-|\cM'|=1$.

By construction $E(\cM)$ contains all but at most $k-1$ edges from each of $P_1$ and $P_2$, so
\[|E(\cM)\setminus E(\cM')|=|E(P_1)\cup E(P_2)\setminus E(P')|\leq 2(k-1)=2(k-1)\mu,\]

Furthermore, $E(\cM')$ contains at most $|E(P')\setminus \lt(E(P_1)\cup E(P_2)\rt)|\leq d+2$ edges that are not contained in $E(\cM)$,
so \[|E(\cM')\setminus E(\cM)|\leq d+2= (d+2)\mu.\]
In conclusion, $\cM'$ is a $(d,k)$-extension of $\cM$,
contradicting the minimality of $\cM$ as a $(d,k)$-extension of itself.
\proofend

\begin{lem}
\label{largematching}
For every $k\geq 1$, in every $n$-vertex graph $G$ satisfying $S(s,g)$ and $L(l)$ for some $s,g$ and $l$ with $sg\geq 4l$,
$s \geq 18$, $\a=\ln s/\ln n$ and $n$ sufficiently large,
for every family of non-trivial vertex-disjoint paths $\cM$ on at most $|V(\cM)|\leq \alpha gs/20$ vertices in $G$
there exists a $(6/\alpha ,k)$-extension $\cM'$ of $\cM$ of size $|\cM'|\leq  \frac{5|V(\cM)|}{2\a ks}+1$.
\end{lem}

\proofstart
Let $\cM'$ be a size-minimum $(6/\alpha,k)$-extension of $\cM$.
(There exists one, since $\cM$ is a $(6/\alpha,k)$-extension of itself.)
Apply Lemma~\ref{inducedexpander} with $D=V(\cM')$ to find the corresponding sets $Z$ and $U$
(meaning that $|Z|\leq 2|V(\cM')|/s$, $U=V(G)\setminus (V(\cM')\cup Z)$ and $G[U]$ satisfies $S(s/2,g/2)$).
Lemma~\ref{inducedexpander} can be applied since
\begin{eqnarray*}
|D|&=&|V(\cM')|=|E(\cM')|+|\cM'|\\
&\leq& |E(\cM)|+(6/\alpha+2)(|\cM|-|\cM'|)+|\cM'|\\
&=& |V(\cM)| + (6/\alpha+1)(|\cM|-|\cM'|) \\
&<& (1+(6/\alpha + 1)/2)|V(\cM)| < 5|V(\cM)|/\alpha\leq gs/4.
\end{eqnarray*}
In the last line we used that the paths of $\cM$ are non-trivial and that $\alpha<1$.

Let now $x$ and $y$ be two vertices from the $k$-ends of two distinct paths $P_1$ and $P_2\in \cM'$, respectively,
and suppose each of them has a neighbor in $U$. Let $a\in U$ be a neighbor of $x$ and let $b\in U$ be a neighbor of $y$.
Since $G[U]$ satisfies both $S(s/2, g/2)$ and $L(l)$ with $l\leq s/2\cdot g/2$,
Observation~\ref{diam} implies that the diameter of $G[U]$ is at most
\[ \diam(G[U])\leq 2\ln n/\ln(s/2)+3\leq 2\frac{\ln n}{\frac{2}{3}\ln s}+3<\frac{6}{\alpha},\]
where the next to last inequality holds since $s\geq 8$.
Hence, there exists an $ab$-path of length at most $\frac{6}{\alpha}$ in $G[U]\subseteq G[V(G)\setminus V(\cM)]$,
contradicting the size-minimality of $\cM'$ by Lemma~\ref{extension}.

Consequently, there can be at most one such path in $\cM'$ that has a vertex in its $k$-end with a neighbor in $U$.
Hence the $k$-ends of at least $|\cM'|-1$ paths have neighbors only in $V(\cM') \cup Z$.
By Lemma~\ref{extension} each such path has length at least $2k-1$
so we found a set of $2k(|\cM'|-1)$ vertices whose neighborhood contains at most
\begin{eqnarray*}
|Z|+|V(\cM')|&<& 2|V(\cM')|/s+|V(\cM')|\leq 10|V(\cM')|/9\\
&<&  \frac{10}{9}(|V(\cM)|    +(6/\alpha+1)|\cM|)\\
& <&5|V(\cM)|/\alpha<gs
\end{eqnarray*}
vertices. Here in the next to last inequality we use the fact that $\cM$ contains no path of length $0$.
Thus, by the small expansion property we obtain
\[|\cM'|-1< \frac{5|V(\cM)|}{2\a ks}.\]
\proofend

The proof of the main theorem is based on the ingenious
rotation-extension technique, developed by P\'osa~\cite{Posa}, and
applied later in a multitude of papers on Hamiltonicity, mostly of
random or pseudorandom graphs (see for example \cite{BFF1},~\cite{FK},~\cite{KS},~\cite{KrSu}).

Let $G$ be a graph and let $P_0=(v_1,v_2,\ldots,v_q)$ be a path in $G$.
If $1 \leq i \leq q-2$ and $(v_q,v_i)$ is an edge of $G$, then there exists a path
$P'=(v_1 v_2\ldots v_i v_q v_{q-1} \ldots v_{i+1})$ in $G$.
$P'$ is called a {\em rotation} of $P_0$ with {\em
fixed endpoint} $v_1$ and {\em pivot} $v_i$. The edge
$(v_i,v_{i+1})$ is called the {\em broken} edge of the rotation. We
say that the segment $v_{i+1} \ldots v_q$ of $P_0$ is {\em reversed} in
$P'$.
In case the new endpoint $v_{i+1}$ has a neighbor $v_j$ such that
$j \notin \{i, i+2\}$, then we can rotate $P'$ further to obtain
more paths of maximum length. We use rotations
together with properties $S$ and $L$ to find a path on the same vertex set as $P_0$ with
large rotation endpoint sets.

The next lemma is a slight strengthening of Claim~2.2 from~\cite{ProfDrSzabo} with a similar proof.
It shows that in any graph having the  small and large expander properties for any path $P_0$ and its endpoint $v_1$
many other endpoints can be created by a small number of rotations with fixed endpoint $v_1$.
In our setting we must also care about not breaking any of the edges from a small ``forbidden'' set $F$.

\begin{lemma}
\label{St}
Let $G=(V,E)$ be a graph on $n$ vertices that satisfies $S(s,g)$ and $L(l)$ with $s\geq 21$, $sg/3>l$, and $l\leq n/24$.
Let $P_0=(v_1,v_2,\ldots,v_q)$ be a path in $G$ and $F\subseteq E(P_0)$ with $|F|\leq s/24-1/2$.
Denote by $B(v_1)\subset
V(P_0)$ the set of all vertices $v\in
V$ for which there is a $v_1v$-path on the vertex set $V(P_0)$ which can be
obtained from $P_0$ by at most $3\frac{\log n}{\log s}$ rotations
with fixed endpoint $v_1$ not breaking any of the edges of $F$.
Then $B(v_1)$ satisfies one of the following properties:
\begin{itemize}
\item
there exists a vertex $v\in B(v_1)$ with a neighbor outside $V(P_0)$, or
\item
$|B(v_1)|\geq n/3$.
\end{itemize}
\end{lemma}

\proofstart
Assume $B(v_1)$ does not have the first property
(i.e., for every $v\in B(v_1)$ it holds that $N(v)\subseteq V(P_0)$).

Let $t_0$ be the smallest integer such that
$\left(\frac{s}{3}\right)^{t_0-2} \geq g$; note that $t_0
\leq  3\frac{\log n}{\log s}$, because $21 \leq s$.

We construct a sequence of sets $S_0, \ldots , S_{t_0}
\subseteq B(v_1) \subseteq V(P_0)\setminus\{v_1\}$ of vertices, such that for
every $0 \leq t \leq t_0$ and every $v\in S_t$, $v$ is the endpoint of a
path which can be obtained from $P_0$ by a sequence of $t$ rotations with fixed endpoint
$v_1$, such that for every $0 \leq i < t$,
the non-$v_1$-endpoint of the path after the $i$th rotation
is contained in $S_i$.
Moreover, $|S_t|=\left(\frac{s}{3}\right)^{t}$ for every $t\leq
t_0-3$, $|S_{t_0-2}|=g$, $|S_{t_0-1}|=l$,
and $|S_{t_0}|\geq n/3$.
Furthermore, for any $1\leq t \leq t_0$, $V(F)\cap S_t=\emptyset$,
and hence no edge from $F$ got broken in any of the rotations.

We construct these sets by induction on $t$. For $t=0$, one can
choose $S_0=\{ v_q\}$ and all requirements are trivially satisfied.

Let now $t$ be an integer with $0< t\leq t_0-2$ and assume that the sets $S_0,
\ldots , S_{t-1}$ with the appropriate properties have already been
constructed. We will now construct $S_t$. Let

\[
T= \{v_i\in N(S_{t-1}) : v_{i-1},v_i,v_{i+1}\not\in
\bigcup_{j=0}^{t-1}S_j\cup V(F)\}
\]
be the set of potential pivots for the $t$th rotation, and notice that $T\subset V(P_0)$ due to our assumption,
since $T\subseteq N(S_{t-1})$ and $S_{t-1}\subseteq B(v_1)$. Assume now
that $v_i\in T$, $y\in S_{t-1}$ and $(v_i,y)\in E$. Then, by the induction hypothesis, a $v_1
y$-path $Q$ can be obtained from $P_0$ by $t-1$ rotations not breaking any edge from $F$ such that
after the $j$th rotation, the non-$v_1$-endpoint is in $S_j$ for
every $0 \leq j \leq t-1$. Each such rotation breaks an edge which is incident with
the new endpoint, obtained in that rotation. Since $v_{i-1}, v_i, v_{i+1}$ are not endpoints
after any of these $t-1$ rotations and also not in $V(F)$, both edges $(v_{i-1},v_i)$ and
$(v_i,v_{i+1})$ of the original path $P_0$ must be unbroken and thus
must be present in $Q\setminus F$.

Hence, rotating $Q$ with pivot $v_i$ will make either $v_{i-1}$ or
$v_{i+1}$ an endpoint (which of the two, depends on whether the unbroken
segment $v_{i-1} v_i v_{i+1}$ is reversed or not after the first
$t-1$ rotations). Assume without loss of generality that the endpoint is $v_{i-1}$. We add $v_{i-1}$
to the set $\hat{S}_{t}$ of new endpoints and say that {\em $v_i$
placed $v_{i-1}$ in $\hat{S}_{t}$}. The only other vertex that can
place $v_{i-1}$ in $\hat{S}_{t}$ is $v_{i-2}$ (if it exists). Thus,

\begin{eqnarray*}
|\hat{S}_{t}| &\ge& \frac{1}{2} |T|\ge
\frac{1}{2} \left(|N(S_{t-1})|-3(1+|S_1|+\ldots+|S_{t-1}|+2|F|\right))\\
&\geq& \frac{s}{2} \left(\frac{s}{3}\right)^{t-1} -
\frac{3}{2}\frac{(s/3)^{t}-1}{s/3-1}-(s/8 - 3/2) \geq
\left(\frac{s}{3}\right)^{t},
\end{eqnarray*}
where the last inequality follows since $s \geq 21$. Clearly we can
delete arbitrary elements of $\hat{S}_{t}$ to obtain $S_{t}$ of size
exactly $\left(\frac{s}{3}\right)^{t}$ if $t\leq t_0-3$ and of size
exactly $g$ if $t=t_0-2$. So the proof of the induction
step is complete and we have constructed the sets $S_0, \ldots,
S_{t_0-2}$.

To construct $S_{t_0-1}$ and $S_{t_0}$ we use the same technique as
above, only the calculations are slightly different.
If $g=1$, then $t_0-1=1$, and analogously to the above calculation we obtain $\hat{S}_{1}$ with $|\hat{S}_{1}|\geq s/3\geq l$.
Otherwise, for $g\geq 2$, since
$|N(S_{t_0-2})|\geq sg$, we have

\begin{eqnarray*}
|\hat{S}_{t_0-1}| &\ge& \frac{1}{2} |T|\ge
\frac{1}{2} \left(|N(S_{t_0-2})|-3(1+|S_1|+\ldots+|S_{t_0-4}|+|S_{t_0-3}|+|S_{t_0-2}|+2|F|\right))\\
&\geq& gs/2 -
\frac{3}{2}\frac{(s/3)^{t_0-2}-1}{s/3-1}-3g/2-(s/8 - 3/2)\\
&\geq&
gs/2 - 2\cdot\left(\frac{s}{3}\right)^{t_0-3} - \frac{3}{2}g-(s/8 - 3/2)\\
&\geq& gs/2 - 2g -\frac{3}{2}g-(s/8 - 3/2)
\geq gs/3> l,
\end{eqnarray*}
where the inequality in the last but one line and the last but one inequality follow since $s \geq 21$ and $g\geq 2$.
We delete arbitrary elements of $\hat{S}_{t_0-1}$ to obtain $S_{t_0-1}$ of size
exactly $l$.

For $S_{t_0}$ the difference in the calculation comes from using the
expansion guaranteed by the property $L$,
rather than the property $S$. That is, we use the fact that
$|N(S_{t_0-1})| \geq n-2l$. Hence, we obtain

\begin{eqnarray*}
|{S}_{t_0}| &\ge& \frac{1}{2} |T|\ge
\frac{1}{2} \left(|N(S_{t_0-1})|-3(1+|S_1|+\ldots+|S_{t_0-2}|+|S_{t_0-1}|+2|F|\right))\\
&\geq& \frac{n}{2}-l -4g-
\frac{3}{2}l-(s/8 - 3/2)\\
&>& \frac{n}{3},
\end{eqnarray*}
where the last inequality follows since $4g\leq sg/3\leq l$, $s\leq sg\leq 3l$ and $l \leq n/24$.

The set $S_{t_0}$ is by construction a subset of $B(v_1)$, concluding the proof of the
lemma. \proofend

Let $H$ be a graph with a spanning path $P=(v_1,\ldots,v_m)$. For
$2\leq i < m$, let us define the auxiliary graph $H^+_i = H^+_{v_i}$ by adding a
vertex and two edges to $H$ as follows: $V(H^+_i)= V(H)\cup \{ w\}$,
$E(H^+_i)=E(H)\cup \{(v_m,w), (v_i,w)\}$. Let $P_i = P_{v_i}$ be the spanning
path of $H^+_i$ which we obtain from the path $P \cup \{(v_m,w)\}$ by
rotating with pivot $v_i$. Note that the endpoints of $P_i$ are
$v_1$ and $v_{i+1}$.

For a vertex $v_i \in V(H)$, let $S^{v_i}$ be the set of those
vertices of $V(P)\setminus \{ v_1\}$, which are endpoints of a
spanning path of $H_i^+$ obtained from $P_i$ by a series of
rotations with fixed endpoint $v_1$.

A vertex $v_i\in V(P)$ is called a {\em bad initial pivot} (or
simply a {\em bad vertex}) if $|S^{v_i}|< \frac{m}{43}$ and is
called a {\em good initial pivot} (or a {\em good vertex})
otherwise. We can rotate $P_i$ and find a large number of endpoints,
provided that $v_i$ is a good initial pivot.

Hefetz et al.~\cite{ProfDrSzabo} showed that $H$ has many good initial pivots
provided that property $L$ is satisfied.

\begin{fact}[\cite{ProfDrSzabo} Lemma 2.3]\label{good}
Let $H$ be a graph satisfying $L(m/43)$ with a spanning path $P=(v_1,\ldots,v_m)$.
Then
\[
|R| \leq 7m/43,
\]
where $R=R(P)\subseteq V(P)$ is the set of bad vertices.
\end{fact}

With these statements in our toolbox, we can prove the following important technical lemma.
It states that we can rotate a path until it can be extended, and still do not break too many of the important edges.

\begin{lem}
\label{path+}
For every sufficiently large $n$ and every $s=s(n)$ with $s\geq 21$,
in every $s$-expander graph $G$ on $n$ vertices every path $P_0$ in $G$ has the following property.
For every pair of sets $F\subset F'\subseteq E(P_0)$ of at most $|F|\leq s/24-1/2$ and $|F'|\leq \frac{n\log s}{9200 \log n}$ edges of $P_0$,
there exists a path $P'$ in $G$ between some $x, y\in V(P_0)$, such that $V(P')=V(P_0)$,
$F\subset E(P')$, $|F'\setminus E(P')|\leq 6\log n/\log s$, and $G$ contains the edge $\{ x,y \} $, or the set
$\{ x,y \} $ has neighbors outside $P'$.
\end{lem}

\proofstart

Assume for the sake of contradiction that the statement is not true.
Let $P_0=(v_1,v_2,\ldots,v_q)$,
and let $A_0=B(v_1)\subset
V(P_0)$ be the set corresponding to $P_0$ and $F$ as in Lemma~\ref{St}, meaning that for every $v\in
B(v_1)$ there is a $v_1v$-path of maximum length which can be
obtained from $P_0$ by at most $t_0=3\frac{\log n}{\log s}$ rotations
with fixed endpoint $v_1$ not breaking any of the edges of $F$.
Clearly, at most $3\frac{\log n}{\log s}$ edges from $F'$ were broken by the rotations,
thus by our assumption every $v\in A_0$ has no neighbors outside $P_0$,
hence by Lemma~\ref{St} we obtain $|A_0|\geq n/3$.
For every $v\in A_0$ fix a $v_1v$-path
$P^{(v)}$ with the above properties and, again using our assumption and Lemma~\ref{St}, construct
sets $B(v)$, $|B(v)| \geq n/3$, of endpoints of paths
with fixed endpoint $v$, obtained from the path $P^{(v)}$ by at most
$t_0$ rotations not breaking any edge from $F$. To summarize, for every $a \in A_0$ and $b \in
B(a)$ there is a path $P(a,b)$ joining $a$ and $b$ on the vertex set $V(P_0)$,
which is obtainable from $P_0$ by at most $\rho:=2t_0 = \frac{6\log
n}{\log s}$ rotations not breaking any of the edges from $F$.
Moreover, this clearly entails $|V(P_0)| \geq n/3$.

We consider $P_0$ to be directed from $v_1$ to $v_q$ and divided
into $2\r$ consecutive vertex disjoint segments $I_1,I_2,\ldots,I_{2\r}$ of length
at least $\rdown{|V(P_0)|/2\r-1}$ each. As every
$P(a,b)$ is obtained from $P_0$ by at most $\r$ rotations, and every
rotation breaks at most one edge of $P_0$, the number of segments of
$P_0$ which also occur as segments of $P(a,b)$, although perhaps reversed,
is at least $\r$. We say that such a segment is {\em unbroken}.
These segments have an absolute orientation given to them by $P_0$,
and another, relative to this one, given to them by $P(a,b)$, which
we consider to be directed from $a$ to $b$. We consider pairs
$\s=(I_{i}, I_{j})$ of unbroken segments of
$P_0$, which occur in this order on $P(a,b)$, where $\s$ also
specifies the relative orientation of each segment.
We call such a
pair $\s$ {\em unbroken}, and say that $P(a,b)$ {\em
contains} $\s$.

\noindent For a given unbroken pair $\s$, we consider the set
$C(\s)$ of ordered pairs $(a,b)$, $a\in A_0, \; b\in B(a)$, such
that $P(a,b)$ contains $\s$.

The total number of unbroken pairs is at most $2^2 (2\r)_2$. Any
path $P(a,b)$ contains at least $\r$ unbroken segments, and thus at
least $\binom{\r}{2}$ unbroken pairs. The average, over
unbroken pairs, of the number of pairs $(a,b)$ such that $P(a,b)$
contains a given unbroken pair is therefore at least
\[\frac{n^2}{9} \cdot \frac{\binom{\r}{2}}{2^2 (2\r)_2} \ge 0.003 n^2.\]
Thus, there is an unbroken pair $\s_0$ and a set $C=
C(\s_0),\,|C|\geq 0.003 n^2$ of pairs $(a,b)$, such that for each
$(a,b)\in C$, the path $P(a,b)$ contains $\s_0$. Let $\hat{A}=\{a
\in A_0 : C$ contains at least $0.003 n/2$ pairs with $a$ as first
element\}. Since $|A_0|, |B(a)| \leq n$, we have $0.003n^2 \leq |C|
\leq |\hat{A}|n + n \cdot \frac{0.003n}{2}$, entailing $|\hat{A}|\geq
0.003 n/2$. For every $a\in \hat{A}$, let $\hat{B}(a)=\{b : (a,b)\in
C\}$. Then, by the definition of $\hat{A}$, for every $a\in\hat{A}$
we have $|\hat{B}(a)|\ge 0.003 n/2$.

For an unbroken pair $\s_0=(I_i,I_j)$,
we divide it
into two segments, $\s_0^1 = (I_i)$
and $\s_0^2 = (I_{j})$, where the parenthesis symbolize that both new
segments maintain the orientation of the segments $I$ in
$\s_0$.
Notice that for every $a\in \hat{A}$ and $b\in \hat{B}(a)$, in the path $P(a,b)$ the segment $\s_0^1$ comes before $\s_0^2$.
For $i=1,2$, let us denote by $|\s_0^i|$ the number of vertices in the
segment $\s_0^i$.
Then for
both segments $\s_0^1$ and $\s_0^2$, we have that
$|\s_0^i|> n/(7\r) $. Let $s_1$ be the first vertex of $\s_0^1$,
$x$ be the last
vertex of $\s_0^1$,
and let $y$ be the first vertex of
$\s_0^2$ and $s_2$ be the last vertex of $\s_0^2$.

We construct a graph $H_1$ with $V(\s_0^1)$ as
vertex set. The edge set of $H_1$ is defined as follows.
First, we add all
edges of $G[V(\s_0^1)]$, except for those that are incident with $s_1$, $x$ or a vertex in $V(F')$.
Further, we add all the edges from $E(\s_0^1)$. Note that all the edges in $H_1$ are also
edges of $G$.
By its construction, $\s_0^1$ is a spanning path in $H_1$ starting at
$s_1$ and ending at $x$.
Let us denote the path reversed to $\s_0^1$ (spanning path in $H_1$,
starting at $x$ and ending at $s_1$) by $P$.
We would like to apply Fact~\ref{good} to $H_1$ with
$m=|V(\s_0^1)|$ and $P$ as the corresponding spanning path.
The condition of the Fact holds since $G$ satisfies property $L(l)$.
Indeed, $l=\frac{n\log s}{3000\log n}$, $m> n/(7\r)=\frac{n\log s}{42\log n}$,
the edges of $H_1$ differ from the edges of $G$ only at $V(F')$ and at the
endpoints of the segment $\s_0^1$, and $|V(F')\cup\{x,s_1\}|\leq \frac{n\log s}{4600 \log n}+2$,
implying
$ l + |V(F')\cup \{x,s_1\}|< \frac{n\log s}{43\cdot 42\log n}<m/43$.
Notice that this are the lines justifying the choices of the integers $3000$ and $9200$ in the bounds for $l$ and $|F'|$.
Hence, $H_1$ satisfies $L(m/43)$, thus by Fact~\ref{good} at least a $\frac{36}{43}$-fraction of the vertices of $H_1$ are good.

For $\s_0^2$ we act similarly: construct a graph $H_2$ from
$\sigma_0^2$ by adding all edges of $G$ with both endpoints in the interior
of $\s_0^2$ but not in $V(F')$ and edges from $E(\s_0^2)$ to $H_2$.
Then $\s_0^2$ forms an oriented spanning path in
$H_2$, starting at $y$ and ending at the last vertex $s_2$ of $\s_0^2$.
Again, due to property $L$,
Lemma~\ref{good} applies here, so
at least a $\frac{36}{43}$-fraction of the vertices of $H_2$ are good.

Recall that $s_1$ is the first vertex of $\s_0^1$.
Since $|\hat{A}| \geq 0.003 n/2 >l +1$ and $H_1$ has at
least $\frac{36}{43}m>l+|V(F')\cup \{x,s_1\}|$
good vertices, there is an edge of $G$ between a vertex
$\hat{a}\in \hat{A}\setminus \{ s_1\}$
and a good vertex $g_1\in V(\s_0^1)\setminus (V(F')\cup \{x,s_1\})$.

Similarly, as $|\hat{B}(\hat{a})|\ge 0.003n/2>l+1$ and there are more than
$l+|V(F')\cup \{y,s_2\}|$ good vertices in $H_2$,
there is an edge from
some $\hat{b}\in \hat{B}(\hat{a})\setminus \{ s_2\}$
to a good vertex $g_2\in V(\s_0^2)\setminus( V(F')\cup \{y,s_2\})$.

Consider the path $P(\hat{a},\hat{b})$ on the vertex set of $P_0$
connecting $\hat{a}$ and $\hat{b}$ and containing $\sigma_0$. The
vertices $x$ and $y$ split this path into three sub-paths: $R_1$
from $\hat{a}$ to $x$, $R_2$ from $y$ to $\hat{b}$ and $R_3$ from
$x$ to $y$. We will rotate $R_1$ with $x$ as a fixed endpoint and
$R_2$ with $y$ as a fixed endpoint, making sure that no edge from $F'$ gets broken.
We will show that the obtained
endpoint sets $V_1$ and $V_2$ are sufficiently large (clearly,
they are disjoint).
Then by property $L$ there will be an edge of $G$ between $V_1$ and $V_2$.
Since we did not touch $R_3$, this edge closes the path into a
cycle, contradicting the assumption from the beginning of the proof.

First we construct the endpoint set $V_1$, the endpoint set $V_2$ can
be constructed analogously.
Recall the notation from Fact~\ref{good}:
Let $H^+_{g_1}$ denote the graph we obtain from $H_1$ by adding the
extra vertex $w$ and the edges $(w,g_1)$ and $(w,s_1)$.
The spanning path of $H^+_{g_1}$
obtained by rotating $P \cup \{ (w,s_1)\}$
with fixed endpoint $x$ at pivot $g_1$ is denoted by $P_{g_1}$.
By the definition of a good vertex, the set $S^{g_1}$
of vertices which are endpoints of a spanning path of $H^+_{g_1}$
that can be obtained from $P_{g_1}$ by a sequence of rotations
with fixed endpoint $x$, has at least $|\s_0^1|/43 >l$ vertices.

We claim that also in $G$, any vertex in $S^{g_1}$ can be obtained as
an endpoint by a sequence of rotations of $R_1$ with fixed endpoint $x$
without breaking any edge from $F'$.
The role of the vertex $w$ will be played by $\hat{a}$ in $G$
(note that we made sure that $\hat{a}\neq s_1$, so $\hat{a}$ is not
contained in $V(\s_0^1)$). Hence, the edge
$(\hat{a},g_1)$ is present in $G$, while we will consider the edge
$(\hat{a},s_1)$ {\em artificial}.

For any endpoint $z\in S^{g_1}$ there is a sequence of pivots, such that
performing the sequence of rotations with
fixed endpoint $x$ at these pivots
results in an $xz$-path spanning $H^+_{g_1}$.
We claim that in $G[V(R_1)]$ it is also possible to perform a
sequence of rotations with the exact same pivot sequence and eventually to
end up in an $xz$-path spanning $V(R_1)$.
When performing these rotations, the subpath of $R_1$
that links $\hat{a}$ to $s_1$
corresponds to the artificial edge $(w,s_1)$ in $H^+_{g_1}$.

Problems in performing these rotations in $G$ could arise if
a rotation is called for where (1)
the pivot is connected to the endpoint of the current spanning path via
an artificial edge of $H^+_{g_1}$: this rotation might not
be possible in $G$ as this edge might not exist in $G$,
or (2)
the broken edge is artificial: after such a rotation in $G$
the endpoint of the new spanning path
might be different from the one we have after
performing the same rotation in $H^+_{g_1}$,
or (3) the broken edge is in $F'$.
However, the construction of $H^+_{g_1}$ ensures that these problems will never occur.
Indeed, in all three cases (1), (2) and (3) the pivot vertex
has an artificial edge or an edge from $F'$ incident with it, while
having degree at least $3$ (as all pivots). However,
both endpoints of an artificial edge and both endpoints of edges from $F'\cap H^+_{g_1}$ have degree $2$ in $H^+_{g_1}$
(for this last assertion we use the fact that $g_1\not\in \{x,s_1\}\cup V(F')$; this is important
as $g_1$ is the first pivot.)

Hence we have ensured that there is indeed a spanning path of $G[V(R_1)]$
from $x$ to every vertex of $V_1=S^{g_1}$ containing all edges from $E(R_1)\cap F'$.

Similarly, since there is an edge from $\hat{b}$ to a good vertex $g_2$
in $H_2$, $g_2\not \in V(F')$, we can rotate $R_2$, starting from this edge to get a set
$V_2=S^{g_2}$ of at least $l$ endpoints not breaking any more edges from $F'$.
In other words we have a spanning path of $G[V(R_2)]$
from $y$ to every vertex of $V_2=S^{g_2}$ containing the edges from $E(R_2)\cap F'$.

As we noted earlier, since by Fact~\ref{good} $|V_1|, |V_2| \geq \frac{m}{43}>l$,
property $L(l)$ ensures that there is an edge between
$V_1$ and $V_2$ in $G$, say $a'b'\in E(G)$ with $a'\in V_1$ and $b'\in V_2$.
This contradicts our assumption, since the rotations we did to obtain $P(a',b')$ from $P_0$ did not break any edges from $F$ and also
all but at most $\r \leq 6\log n/\log s$ edges from $F'$ are on the path $P(a', b')$.
\proofend

We are now able to prove the main lemma, stating that for every matching there is a Hamilton cycle almost covering it.
Notice that we use the same calligraphic letter $\cM$ to denote a matching as for families of paths,
since we are going to apply extension on the matching and hence we see it as a family of paths of length $1$ each.

\begin{lem}
\label{mainlemma}
For every constant $\alpha\in (0,1]$ and for every sufficiently large $n$ the following holds.
Let $G$ be an $n^\alpha$-expander graph on $n$ vertices.
For every matching $\cM$ in $G$ of size at most $|\cM|\leq \alpha^3n/9200$
there exists a Hamilton cycle $C$ in $G$ with
\[|E(\cM)\setminus E(C)|\leq \left\lfloor \frac{1036 |\cM|}{\alpha^3n^{\alpha/2}}\right\rfloor\].
\end{lem}

\proofstart
First we proceed inductively to construct a single path via $(d,k)$-extensions that contains
most of the matching edges.

Using Lemma~\ref{largematching} with $s=n^\alpha$, $g=4\alpha n^{1-\alpha}$, $l=  \alpha n/3000  $, $k=1$ and $d=6/\alpha$,
we find a $(d,1)$-expansion $\cM_2$ of $\cM_1=\cM$ of size at most $\left\lfloor \frac{5|\cM|}{\alpha n^\alpha}+1\right\rfloor$
containing all edges of $\cM$.

For $i\geq 2$, given a family of vertex-disjoint non-trivial paths $\cM_i$ of size
\[2\leq |\cM_i| \leq \left\lfloor \frac{45 |\cM|}{2\alpha^2n^{i\alpha/2}}+1\right\rfloor \]
on at most
\[|V(\cM_i)|\leq (d+3)|\cM| \] vertices
containing all but at most
\[(i-2)\frac{45 |\cM|}{\alpha^2n^{\alpha/2}}\]
edges from $\cM$,
we construct a size-minimum $(d, n^{(i-1)\alpha/2})$-extension $\cM_{i+1}$ of $\cM_i$.
Then $\cM_{i+1}$ satisfies the above properties by construction:
it contains all but at most
\begin{eqnarray*}
|E(\cM)\setminus E(\cM_{i+1})|&\leq & |E(\cM_{i})\setminus E(\cM_{i+1})|+|E(\cM)\setminus E(\cM_{i})| \\
&\leq & 2n^{(i-1)\alpha/2}(|\cM_i|-1)+(i-2)\frac{45 |\cM|}{\alpha^2n^{\alpha/2}}\\
&\leq & (i-1)\frac{45 |\cM|}{\alpha^2n^{\alpha/2}}
\end{eqnarray*}
edges from $\cM$.
Furthermore, since $\cM_{i+1}$ was constructed from $\cM$ by a series of $(d,k)$-extensions with varying $k$ but fixed $d$,
at most $d+1$ vertices were added between any two paths to concatenate them,
thus $\cM_{i+1}$
has at most $(d+3)|\cM|$
vertices.

Finally,
by Lemma~\ref{largematching} $\cM_{i+1}$ has size at most
\[  |\cM_{i+1}|\leq \left\lfloor \frac{5|V(\cM_i)|}{2\alpha n^{(i-1)\alpha/2}n^\a}+1\right\rfloor
 \leq \left\lfloor \frac{5(d+3) |\cM|}{2\a n^{(i+1)\alpha/2}}+1\right\rfloor
 \leq \left\lfloor \frac{45 |\cM|}{2\alpha^2n^{(i+1)\alpha/2}}+1\right\rfloor  . \]

The family $\cM_{last}$, where $last\leq 2/\alpha+1<3/\alpha$ is the index we stop the induction with, contains only one path $P$.
Let us apply Lemma~\ref{inducedexpander} with $D=V(P)$. This is possible, since $P$ contains at most
$(d+3)|\cM|<\alpha n =gs/4$ vertices.
We obtain the corresponding sets $Z$ and $U=V\setminus (D\cup Z)$ and conclude that  the induced graph $G[U]$
satisfies the small expander property $S(s/2, g/2)$.
Theorem 2.5 from~\cite{ProfDrSzabo} states that for every choice of the expansion parameter $r$
with $12\leq r \leq \sqrt{n}$,
every $n$-vertex graph $G$ satisfying
$S\lt(r,\frac{n\ln r}{r\ln n}\rt)$ and $L\lt(\frac{n\ln r}{1035 \ln n}\rt)$ is Hamiltonian.
Hence, applying this statement to $G[U]$
with $s=n^\alpha$, $g=4\alpha n/n^\alpha$ and $r=n^\alpha/2$,
we see that $G[U]$ is Hamiltonian.

Furthermore, by Lemma~\ref{inducedexpander} we know that
\begin{align}
|Z|\leq 2|V(P)|/n^\a. \label{|z|}
\end{align}
Using the small expander property of $G$,
we obtain an edge between a vertex $x$ in the $\left\lceil |\cM|/n^{\alpha/2}\right\rceil$-end of $P$
and a vertex $y$ in $U$ in the following way.
Take a subset of size $\min\{2\left\lceil |\cM|/n^{\alpha/2}\right\rceil, 4\a n^{1-\a}\}$ of the $\left\lceil |\cM|/n^{\alpha/2}\right\rceil$-end of $P$.
Its neighborhood has size at least
$\min\{2\left\lceil |\cM|/n^{\alpha/2}\right\rceil n^\a, 4\a n- |V(P)|\}>|V(P)\cup Z|$.

The vertex $x$ breaks the path $P$ into two subpaths, one of which contains all but at most $|\cM|/n^{\alpha/2}$ edges from $P$.
Connecting this path via the edge $xy$ with a Hamilton path in $G[U]$,
we create a path $R$ containing all but at most
\begin{align}
|\cM\setminus E(R)|&\leq
\left\lfloor\frac{3}{\a}\cdot \frac{ 45 |\cM|}{\alpha^2n^{\alpha/2}}\right\rfloor+\left\lfloor|\cM|/n^{\alpha/2}\right\rfloor
\leq \left\lfloor\frac{136 |\cM|}{\alpha^3n^{\alpha/2}}\right\rfloor \label{m-p}
\end{align}
edges from $\cM$
and all vertices from $U$.
Notice that from~(\ref{|z|}) we have that $R$ is missing only
\begin{align}
n-|E(R)|\leq |Z|+|\cM|/n^{\alpha/2} +1< 2|\cM|/n^{\alpha/2}+1.\label{n-r}
\end{align}
vertices to be Hamiltonian.

We aim to use Lemma~\ref{path+} to rotate/extend $R$ into a Hamilton cycle without losing
many edges of $\cM$ that are already on it.

For this we set $P_0=R$ and $F'=F=\cM$, in case $|\cM|<\alpha^3n^{\alpha/2}/136$
(and thus $\cM$ is contained in $R$ by~(\ref{m-p})).
Otherwise, if $|\cM|\geq\alpha^3n^{\alpha/2}/136$, let $F=\emptyset$ and $F'=\cM\cap E(R)$.

We will now use Lemma~\ref{path+} iteratively, in each step rotating/extending our current path with an edge until
it is spanning and then closing it into a Hamilton cycle.
Notice that Lemma~\ref{path+} can be applied throughout the process
since $|F'|\leq |\cM|<\alpha n/9200$ and $|F|=o(n^\alpha)$.

Consider the $x,y$-path $P'$ arising from the application of  Lemma~\ref{path+} to $F$, $F'$ and $P_0$.
If one of the two vertices $x$ or $y$ has neighbors outside $P'$,
we can extend $P'$ with one more edge to obtain a longer path $\hat{P}$ containing $P'$.
We update for our iteration $P_0:=\hat{P}$ and $F':=\cM \cap E(\hat{P})$ and start the iteration step again.
Notice that in this step, the size of $F'$ decreased by at most $6/\alpha$.
If $x$ and $y$ have no neighbors outside $P'$, then there is a cycle $C$ containing $P'$.
If $C$ is Hamilton, we stop the procedure since this is what we are aiming at.
Otherwise, by the connectivity of $G$ (guaranteed by properties $S$ and $L$ and stated implicitly in Observation~\ref{diam}),
there is a vertex $w\in V(C)$ with a neighbor outside $C$, say $aw\in E(G)$, $a\in V(G)\setminus V(C)$.
Notice that only one of the edges incident with $w$ in $C$ can be in $F'$, since $F'\subset \cM$ is a matching.
Removing an edge incident with $w$ in $C$ which is not in $F'$ and adding the edge $aw$,
we obtain a path $\hat{P}$ of length $|\hat{P}|\geq|P'|+1$ containing all edges from $F'\cap E(P')$.
We update for our iteration $P_0:=\hat{P}$ and $F':=\cM\cap E(\hat{P}) $ and start the iteration step again.
Notice that again, in this step, the size of $F'$ decreased by at most $6/\alpha$.

Using~(\ref{n-r}), we see that after at most $n-|R|\leq 2|\cM|/n^{\alpha/2}+1$ steps the iteration ends and we obtain a Hamilton cycle $C$.
If $|\cM|<\alpha^3n^{\alpha/2}/136$, then $C$ contains all edges from $\cM$.  
Otherwise, 
$C$ contains all but at most
\begin{align*}
|\cM\setminus E(C)|&\leq |E(R)\setminus E(C)|+|\cM\setminus E(R)|
\leq \lt(2|\cM|/n^{\alpha/2}+1\rt)\cdot \frac{6}{\alpha}+\frac{136 |\cM|}{\alpha^3n^{\alpha/2}}\\
&=\frac{6}{\a}+\frac{12|\cM|}{\a n^{\a/2}}+\frac{136 |\cM|}{\alpha^3n^{\alpha/2}}
< \frac{7}{\a}\cdot \frac{148 |\cM|}{\alpha^3n^{\alpha/2}}
\end{align*}
edges from $\cM$, completing the proof of the lemma.

\proofend

The following corollary condenses all the previous technical work.
It states that every matching of an $n^\alpha$-expander graph can be covered with a constant-size collection of Hamilton cycles.

\begin{cor}
\label{oldteo}
For every constant $\alpha, 0<\alpha \leq 1$, there is an $n_0$, such that
for every $n\geq n_0$ the following holds. In every $n^\alpha$-expander graph $G$ on $n$ vertices,
for every matching $M$ of $G$ there exist at most $14000/\alpha^4$ Hamilton cycles such that
$M$ is contained in their union.
\end{cor}
\proofstart
We start by splitting $M$ into at most $\lt\lceil 4600/\a^3\rt\rceil<4601/\a^3$ matchings of size at most $\alpha^3n/9200$ each.
For every such matching $M_1$ we set $i=1$ and perform the following iterative procedure:
\begin{itemize}
\item
take a Hamilton cycle $C_i$ covering as many of the edges of $M_i$ as possible.
If $M_i\subset E(C_i)$, then we found our covering and finish the procedure.
\item
otherwise, we set $M_{i+1}:=M_i\setminus E(C_i)$ and remark that by Lemma~\ref{mainlemma}
\begin{align}
|M_{i+1}|\leq \frac{1036 |M_i|}{\alpha^3n^{\alpha/2}}\leq \lt(\frac{1036}{\alpha^3n^{\alpha/2}}\rt)^{i+1}|M_1|.\label{count}
\end{align}
Update $i:=i+1$ and start again from the first iteration step.
\end{itemize}
From~(\ref{count}) we see that after at most $\lfloor2/\a +1\rfloor<3/\a$ iteration steps,
we get a collection of at most $3/\a$ Hamilton cycles $C_1, C_2,\ldots$ covering $M_1$.
Hence, there exist a collection of at most $\frac{3}{\a}\cdot \frac{4601}{\a^3}$ Hamilton cycles covering $M$,
implying the statement of the lemma.
\proofend

We are now able to prove Theorem~\ref{teo}.

\proofstart

We start by taking $h$ disjoint Hamilton cycles into our covering.
Removing the union of these cycles from $G$, we are left with a graph $H$
of maximum degree exactly $\Delta(H)=\Delta(G)-2h$.
Using at most $2\Delta(H)$ colors we color the edges of $H$ greedily, partitioning them into at most $2\Delta(H)$ matchings.
By Corollary~\ref{oldteo} for every of these matchings there exist $14000/\a ^4$ Hamilton cycles covering it,
completing the proof of the theorem.

\proofend

\section{Proof of Theorem~\ref{corrandom}}\label{pf2}

In this section we derive the proof of Theorem~\ref{corrandom}  from Theorem~\ref{teo}
by checking that $G(n,p)$ is an $s$-expander a.a.s. with the appropriate choice of $p$ and $s$.
Notice that this choice of $s$ is clearly not optimal, but suffices for our purposes.

\begin{lem}
\label{expander}
For every constant $\alpha$ with $0< \alpha <1$ and every function $p=p(n)\geq n^{\alpha-1}$,
$G(n,p)$ is a $\sqrt[5]{np}$-expander a.a.s.
\end{lem}

\proofstart
Let $s=\sqrt[5]{np}$.

First we prove that $G(n,p)$ has the small expander property $S(s,  \frac{4n\ln s}{s\ln n})$.

Let $A\subset V(G(n,p))$ be an arbitrary subset of size at most $|A|\leq \frac{4n\ln s}{s\ln n} \leq \frac{4n}{5s}$.
The random variable $|N(A)|$ is the sum of the $n-|A|$ characteristic variables
of the events $v\in N(A)$ for $v \in V\backslash A$.
Hence for the expectation we obtain
\begin{align*}
\Exp[|N(A)|] & =\sum_{v\in V\setminus A} \Pr[ v\in N(A)]=(n-|A|)\left(1-(1-p)^{|A|}\right)>(n-|A|)\frac{|A|p}{1+|A|p}\\
&\geq(1+o(1))n\frac{|A|p}{1+|A|p} \geq  (1+o(1))|A|s\frac{np/s}{1+\frac{4}{5}np/s} = |A|s\left(\frac{5}{4} +o(1)\right).
\end{align*}
Here we first used the simple fact that $(1-p)^{|A|}<\frac{1}{1+|A|p}$, then a couple of times  that $|A|\leq \frac{4n}{5s}$.
Since the elementary events $v\in N(A)$ that make up $|N(A)|$ are mutually independent the Chernoff bound can be applied
to estimate the probability that $A$ is not expanding. We use the above estimate on $\Exp[|N(A)|]$ several times.
\begin{align*}
\Pr\left[|N(A)|  \leq  s|A|\right]
&< \exp\left[-\frac{(\Exp[|N(A)|]-|A|s)^2}{2\Exp[|N(A)|]}\right]
< \exp \left( -\frac{\lt(\lt(\frac{1}{5}+o(1)\rt)\Exp[|N(A)|]\rt)^2}{2\Exp[|N(A)|]}\right)\\
&= \exp \left( -\left(\frac{1}{50} +o(1)\right)\Exp[|N(A)|]\right)
< \exp \left( -\left(\frac{1}{40} +o(1)\right)|A|s\right).
\end{align*}

By the union bound the probability that $G(n,p)$ does not satisfy property $S\lt(s,  \frac{4n\ln s}{s\ln n}\rt)$ is bounded by
\begin{align*}
\Pr\left[\exists A\subset V, |A|\leq \frac{4n}{5s}:~|N(A)|  \leq  s|A|\right]
&< \sum_{a=1}^{n}\binom{n}{a}\exp\left[-\lt(\frac{1}{40}+o(1)\rt)a s\right]\\
&< \sum_{a=1}^{\infty}\left(n\exp\left[-0.01 n^{\alpha/5}\right]\right)^{a}=o(1).
\end{align*}

To complete the proof we show that $G(n,p)$ has the large expander property $L\lt(\frac{\a n}{15000 \ln n}\rt)$.

Let $A,B\subseteq V(G(n,p))$ be fixed subsets of size $|A|,|B|\geq\frac{\a n}{15000 \ln n}$ with $A\cap B= \emptyset$.
Then we have that

\[\Pr[\mbox{there are no edges between } A \mbox{ and }B]= (1-p)^{|A||B|}<\exp\lt(-\frac{\a^2n^2p}{15000^2 \ln^2 n}\rt).\]

Using union bound over all pairs of such disjoint sets $A,B\subseteq V(G(n,p))$, we get the desired probability
\[\Pr\lt[G(n,p)\mbox{ satisfies } L\lt(\frac{\a n}{15000 \ln n}\rt)\rt]
\geq 1-4^n \exp\lt(-\frac{\a^2n^2p}{15000^2 \ln^2 n}\rt)=1-o(1),\]
proving the lemma.
\proofend

We are now able to prove Theorem~\ref{corrandom} using Theorem~\ref{teo}.

\proofstart Let $p$ be in the range of the theorem.
By Lemma~\ref{expander} $G(n,p)$ is an $n^{\alpha /5}$-expander a.a.s.
We know from~\cite{KKO1} that there exists a packing of $(1-o(1))np/2$ Hamilton cycles into $G(n,p)$.
Finally, the maximum degree a.a.s. satisfies $\Delta(G(n,p))=(1+o(1))np$.
Hence, by Theorem~\ref{teo}, we obtain a covering of $G(n,p)$ by
$(1-o(1))np/2+28000((1+o(1))np-2(1-o(1))np/2)/(\alpha/5)^4=(1+o(1))np/2$ Hamilton cycles, finishing the proof of the theorem.
\proofend

\section{Concluding remarks and open questions}\label{concl}



In this paper we verified that the size of a largest Hamilton cycle
packing and the size of a smallest Hamilton cycle covering are
asymptotically equal a.a.s. in the random graph $G(n,p)$, provided
$p\geq n^{\alpha -1}$ for an arbitrary constant $\alpha >0$. Our
result calls for at least two natural directions of possible
improvement.

First of all, the only  explanation why the lower bound on the edge probability needs to be at least $n^{\alpha -1}$
and the corresponding expansion factor $s$ needs to be at least $n^{\a/5}$ is that
for lower values of $p$ and $s$ most of our technical arguments would break down.
We think these bounds on $p$ and $s$ are only artifacts of our proof.
Since the minimum and maximum degrees of the random graph have to be asymptotically equal
in order for the minimum Hamilton covering and the maximum Hamilton packing to be asymptotically
of the same size, we need to assume that $p=\omega(\ln n)/n$.
We strongly believe though that Theorem~\ref{corrandom} 
holds already whenever $p=\omega(\ln n)/n$. 
\begin{conj}
For any $p=\omega (\ln n)/n$ the random graph $G(n,p)$ admits a
covering of its edges with at most $(1+o(1))np/2$ Hamilton cycles
a.a.s.
\end{conj}
Even though~\cite{KKO1} provides the corresponding packing result for this range,
we were not able to extend our techniques to prove e.g. the analog of Lemma~\ref{mainlemma} for $\alpha=o(1)$,
not to mention for  $\a=(\ln\ln n+\omega(1))/\ln n$.

Another direction in which Theorem~\ref{corrandom} could be
tightened is to make the statement exact instead of approximate. The
trivial lower bound on the size of a Hamilton covering in terms of
the maximum degree
is $\lceil \Delta (G)\rceil$. In what range of $p$ will this be tight.

\begin{question}
\label{q1}
In what range of $p$ does there exist a Hamilton covering of $G\sim G(n,p)$ of size $\lceil\Delta(G)/2\rceil$ a.a.s.?
\end{question}
Recall that the analogous precise statement in terms of the minimum
degree {\em is} true for Hamilton packings \cite{KKO2}, \cite{KrSa}.
To have a positive answer for the question, we clearly need to be
above the Hamiltonicity threshold, but it is plausible that the
statement is true immediately after that.
%

The question of covering the edges of a graph by Hamilton cycles can
also be considered in the {\em pseudorandom} setup. A graph $G$ is
called an $(n, d, \lambda)$-graph if it is $d$-regular on $n$
vertices and the second largest absolute value of its eigenvalues is
$\lambda$.
The concept of $(n,d,\lambda)$-graphs is a common way to formally
express pseudorandomness, as $(n,d,\lambda)$-graphs with
$\lambda=o(d)$ behave in many ways as random graphs are expected to
do. (See, e.g., ~\cite{KrSu2} for a general discussion on
pseudorandom graphs and $(n,d,\lambda)$-graphs.)
Theorem 2 from~\cite{FK1} implies that for $(n,d,\lambda)$-graph with $d=\Theta(n)$ and $\lambda=o(d)$
there exists a Hamilton packing of size $d/2-3\sqrt{\lambda n}=d/2-o(d)$.
Then using the result of Tanner~\cite{tanner} that in every $(n,d,\lambda$-graph $G$ for every subset $X$ of $V(G)$,
\[|N(X)|\geq \frac{d^2|X|}{ \lambda^2+(d^2-\lambda^2)|X|/n},\]
one can see by a simple case analysis that whenever there exists a
constant $0<\a\leq 1/5$ such that $d=\Theta(n)$ and
$\lambda=O\lt(n^{1-\a}\rt)$, then every $(n,d,\lambda)$-graph $G$
with sufficiently large $n$ is an $n^\a$-expander.
Hence, Theorem~\ref{teo} implies that any such $(n,d,\lambda)$-graph
has a Hamilton covering of size $d/2+o(d)$, which is of course
asymptotically best possible. It would be interesting to decide
whether a similar statement holds for sparser pseudorandom graphs,
maybe as sparse as $d=n^\epsilon$, for arbitrarily small $\epsilon
>0$.

A Hamilton cycle is a particular spanning structure of the complete
graph, which can be used to decompose its edges. A further group of
problems related to our result is to determine the typical sizes of
a largest packing and of a smallest covering of various other
spanning structures in the random graph. Here often the
corresponding decomposition result for the complete graph is not
known or is just conjectured. Still asymptotic packing and covering
results would be of interest, for example for trees of bounded
maximum degree.

\end{document}